\documentclass[11pt]{amsart}
\usepackage{amssymb,latexsym}
\usepackage{pstcol,pstricks,color}

 \numberwithin{equation}{section}

\newtheorem{theorem}{Theorem}[section]
\newtheorem{proposition}[theorem]{Proposition}
\newtheorem{corollary}[theorem]{Corollary}

\newcommand{\csum}{{\rm csum}}
\newcommand{\sumlim}{\sum\limits}

\newcommand{\een}{\end{enumerate}}
\newcommand{\blem}{\begin{lem}}
\newcommand{\elem}{\end{lem}}
\newcommand{\bcl}{\begin{clm}}
\newcommand{\ecl}{\end{clm}}
\newcommand{\bthm}{\begin{thm}}
\newcommand{\ethm}{\end{thm}}
\newcommand{\bpr}{\begin{prop}}
\newcommand{\epr}{\end{prop}}
\newcommand{\bco}{\begin{cor}}
\newcommand{\eco}{\end{cor}}
\newcommand{\bcon}{\begin{conj}}
\newcommand{\econ}{\end{conj}}
\newcommand{\bde}{\begin{defn}}
\newcommand{\ede}{\end{defn}}
\newcommand{\bex}{\begin{exa}}
\newcommand{\eexa}{\end{exa}}
\newcommand{\bobs}{\begin{obs}}
\newcommand{\eobs}{\end{obs}}
\newcommand{\bexe}{\begin{exe}}
\newcommand{\eexe}{\end{exe}}

\newcommand{\grn}{G_{r,n}}

\def\qed{\hfill $\Box$}
\def\pf{\noindent {\it Proof.} }

\title{Excedance numbers for permutations in complex reflection groups}

\begin{document}
\maketitle
\begin{center}
Toufik Mansour$^\dag$ and Yidong Sun$^\ddag$

$^\dag$Department of Mathematics, University of Haifa, 31905
Haifa, Israel\\
$^\ddag$Department of Mathematics, Dalian Maritime University, 116026 Dalian, P.R. China\\[5pt]

{\it $^\dag$toufik@math.haifa.ac.il, $^\ddag$sydmath@yahoo.com.cn}
\end{center}\vskip0.5cm

\subsection*{Abstract}
Recently, Bagno, Garber and Mansour \cite{BGM} studied a kind of
excedance number on the complex reflection groups and computed its
multidistribution with the number of fixed points on the set of
involutions in these groups. In this note, we consider the similar
problems in more general cases and make a correction of one result
obtained by them.

{\bf Keywords}: Complex reflection group, excedance, colored
permutations

\noindent {\sc 2000 Mathematics Subject Classification}: Primary
05A05, 05A15

\section{Introduction}
It is well known, there is a single infinite family of groups
$G_{r,s,n}$ and exactly $34$ other "exceptional" complex reflection
groups. The infinite family $G_{r,s,n}$, where $r,s,n$ are positive
integers with $s\mid r$, consists of the groups of $n\times n$
matrices such that
\begin{itemize}
\item the entries are either 0 or $r^{\rm th}$ roots of unity;
\item  there is exactly one nonzero entry in each row and each column;
\item  the $(r/s)^{{\rm th}}$ power of the product of the nonzero
entries is 1,
\end{itemize}
 where $r,s,n$ are positive integers with $s\mid r$.

The classical Weyl groups appear as special cases: for $r=s=1$ we
have the symmetric group $G_{1,1,n}=S_n$, for $r=2s=2$ we have the
hyperoctahedral group $G_{2,1,n}=B_n$, and for $r=s=2$ we have the
group of even-signed permutations $G_{2,2,n}=D_n$.

We say that a permutation $\pi \in G_{r,s,n}$ is {\it an involution}
if $\pi^2=1$. More generally, we define
$\mathcal{G}_{r,s,n}^{m}=\{\sigma\in G_{r,s,n}|\sigma^m=1\}$.
Recently, Bagno, Garber and Mansour \cite{BGM} studied an excedance
number on the complex reflection groups (see~\cite{SS}) and computed
the number of involutions having specific numbers of fixed points
and excedances. In this note, we consider the similar problems on
the set $\mathcal{G}_{r,s,n}^{m}$.

This paper is organized as follows. In Section \ref{pre}, we recall
some properties of $G_{r,s,n}$ and define some parameters on
$\grn=G_{r,1,n}$ and hence also on $G_{r,s,n}$. In Section
\ref{mainresults} we present our main results and compute the
corresponding recurrences together with explicit formulas.

%=====================================
\section{Preliminaries}\label{pre}
Let $r$ and $n$ be any two positive integers. {\it The group of
colored permutations of $n$ digits with $r$ colors} is the wreath
product $\grn=\mathbb{Z}_r \wr S_n=\mathbb{Z}_r^n \rtimes S_n$
consisting of all the pairs $(z,\tau)$ where $z\in\mathbb{Z}_r^n$
and $\tau \in S_n$. Let $\tau,\tau'\in S_n$,
$z=(z_1,...,z_n)\in\mathbb{Z}_r^n$ and
$z'=(z'_1,...,z'_n)\in\mathbb{Z}_r^n$, the multiplication in $\grn$
is defined by
\begin{eqnarray*}
(z,\tau) \cdot
(z',\tau')=((z_1+z'_{\tau^{-1}(1)},...,z_n+z'_{\tau^{-1}(n)}),\tau
\circ \tau'),
\end{eqnarray*}
where $+$ is taken modulo $r$.

We use some conventions along this paper. For an element
$\sigma=(z,\tau) \in \grn$ with $z=(z_1,...,z_n)$ we write
$z_i(\sigma)=z_i$. For $\sigma=(z,\tau)$, we denote
$|\sigma|=(0,\tau), (0 \in \mathbb{Z}_r^n)$.

A much more natural way to present $\grn$ is the following: Consider
the alphabet
$\Sigma=\{1^{[0]},\dots,n^{[0]},1^{[1]},\dots,n^{[1]},\dots,
1^{[r-1]},\dots,n^{[r-1]} \}$ as the set $[n]=\{1,\dots,n\}$ colored
by the colors $0,\dots,r-1$. Then, an element of $\grn$ is a {\it
colored permutation}, i.e. a bijection $\sigma: \Sigma \rightarrow
\Sigma$ such that if $\sigma(i)=k^{[t]}$ then
$\sigma(i^{[j]})=k^{[t+j]}$ where $0\leq j \leq r-1$ and the
addition is taken modulo $r$. Occasionally, we write $j$ bars over a
digit $i$ instead of $i^{[j]}$. For example, an element
$(z,\tau)=((1,2,1,2),(3,1,2,4)) \in G_{3,4}$ will be written as
$(\bar{3} \bar{\bar{1}} \bar{2} \bar{\bar{4}})$.

For each $s|r$ we define the {\it complex reflection group}:
$$G_{r,s,n}:=\{\sigma \in G_{r,n} \mid \csum(\sigma)\equiv 0  \; {\rm
mod} \; s\},$$ where ${\rm csum}(\sigma) = \sumlim_{i=1}^n
z_i(\sigma)$.

One can define the following well-known statistics on $S_n$. For any
permutation $\sigma \in S_n$, $i \in [n]$ is {\it an excedance of
$\sigma$} if and only if $\sigma(i)>i$. We denote the number of
excedances by ${\rm exc}(\sigma)$. Another natural statistic on
$S_n$ is the number of fixed points, denoted by ${\rm fix}(\sigma)$.
We can similarly define some statistics on $\grn$. The complex
reflection group $G_{r,s,n}$ inherits all of them. Given any ordered
alphabet $\Sigma'$, we recall the definition of the {\it excedance
set} of a permutation $\sigma$ on $\Sigma'$:
$${\rm Exc}(\sigma)=\{i \in \Sigma' \mid \sigma(i)>i\},$$ and the {\it excedance
number} is defined to be ${\rm exc}(\sigma)=|{ \rm Exc}(\sigma)|$.

We define the color order on the set
$\Sigma=\{1,\dots,n,\bar{1},\dots,\bar{n},\dots,
1^{[r-1]},\dots,n^{[r-1]} \}$ for $0\leq j<i< r$ by $1^{[i]}<2^{[i]}
< \cdots < n^{[i]} < 1^{[j]} < 2^{[j]} < \cdots < n^{[j]}$. We note
that there are some other possible ways of defining orders on
$\Sigma$, some of them lead to other versions of the excedance
number, see for example \cite{BG}. For example, given the color
order $\bar{\bar{1}} < \bar{\bar{2}} <\bar{\bar{3}} < \bar{1} <
\bar{2} <\bar{3} < 1 < 2 < 3$, we write
$\sigma=(2\bar{1}\bar{\bar{3}}) \in G_{3,3}$ in an extended form
$$(\star)\hskip2cm\begin{pmatrix} \bar{\bar{1}} & \bar{\bar{2}} & \bar{\bar{3}} &
\bar{1} & \bar{2}& \bar{3} & 1 & 2 & 3\\
\bar{\bar{2}} & 1 & \bar{3} & \bar{2} & \bar{\bar{1}}  &  3 & 2 &
\bar{1} & \bar{\bar{3}}
\end{pmatrix}$$
which implies that  ${\rm
Exc}(\sigma)=\{\bar{\bar{1}},\bar{\bar{2}},\bar{\bar{3}},\bar{1},\bar{3},1\}$
and ${\rm exc}(\sigma)=6$.

Define ${\rm Exc}_A(\sigma) = \{ i \in [n-1] \ | \ \sigma(i) > i
\}$, where the comparison is with respect to the color order, and
denote ${\rm exc}_A(\sigma) = |{\rm Exc}_A(\sigma)|$. For instance,
if $\sigma=(\bar{1}\bar{\bar{3}}2\bar{\bar{4}}) \in G_{3,4}$, then
${\rm csum}(\sigma)=5$, ${\rm Exc_A}(\sigma)=\{3\}$ and hence ${\rm
exc}_A(\sigma)=1$.

Now we can define the colored excedance number for $\grn$ by ${\rm
exc}^{{\rm Clr}}(\sigma)=r \cdot {\rm exc}_A(\sigma)+{\rm
csum}(\sigma)$. Let $\Sigma$ ordered by the color order, then we can
state that ${\rm exc}(\sigma)={\rm exc}^{{\rm Clr}}(\sigma)$
obtained by Bagno and Garber \cite{BG} for any $\sigma \in \grn$.

For $\sigma=(z,\tau)\in \grn$, $|\sigma|$ is the permutation of
$[n]$ satisfying $|\sigma|(i)=\tau(i)$. We say that $i \in [n]$ is
an {\it absolute fixed point} of $\sigma \in \grn$ if
$|\sigma|(i)=i$. We denote the number of absolute fixed point of
$\sigma \in \grn$ by ${\rm fix(\sigma)}$.

\section{Main results and proofs}\label{mainresults}
Our main result can be formulated as follows.  Recall that
$\mathcal{G}_{r,s,n}^{m}=\{\sigma\in G_{r,s,n}|\sigma^m=1\}$, define
\begin{eqnarray*}
H_{r,s,n}^{(m)}(u,v,w)=\sum_{\sigma\in
\mathcal{G}_{r,s,n}^{m}}u^{{\rm fix}(\sigma)}v^{{\rm
exc}_{A}(\sigma)}w^{{\rm csum}(\sigma)},
\end{eqnarray*}
and
\begin{eqnarray*}
\mathcal{H}_{r,s}^{(m)}(x;u,v,w)=\sum_{n \geq
0}{H_{r,s,n}^{(m)}(u,v,w)\frac{x^n}{n!}}=\sum_{n\geq0}\sum_{\sigma\in
\mathcal{G}_{r,s,n}^{m}}\left( u^{{\rm fix}(\sigma)}v^{{\rm
exc}_A(\sigma)}w^{{\rm csum}(\sigma)}\right)\frac{x^n}{n!}.
\end{eqnarray*}

\begin{theorem}\label{maintheo}
For positive integer $r,m\geq 1$, there holds
$$\begin{array}{l}
\mathcal{H}_{r,1}^{(m)}(x;u,v,w)\\
\qquad=\exp\left\{\sum\limits_{\{t|0\leq t<r, r|tm\}}xuw^t
+\sum\limits_{d|m,d\geq
2}\frac{x^d}{d!}\sum\limits_{k=1}^{d-1}A_{d-1,k}\sum\limits_{i=0}^k\binom{k}{i}v^{k-i}\sum\limits_{r|\frac{tm}{d}}
U_{d-k,t}^{(i)}w^{t}\right\},\end{array}$$ where $A_{d-1,k}$ is the
Eulerian number, that is the number of permutations on $[d-1]$ with
$k-1$ excedances, $U_{d-k,t}^{(i)}$ is the coefficient of $x^t$ in
$(x+x^2+\cdots+x^{r-1})^i(1+x+\cdots+x^{r-1})^{d-k}$, namely,
\begin{eqnarray*}
U_{d-k,t}^{(i)}=\sum_{j=0}^i(-1)^{i-j}\binom{i}{j}\sum_{{\ell}\geq
0}(-1)^{{\ell}}\binom{d+j-k}{{\ell}}\binom{d+j+t-k-{\ell}r-1}{t-{\ell}r}.
\end{eqnarray*}
\end{theorem}
\pf For any $\pi \in \mathcal{G}_{r,1,n}^{m}$, the length of each
cycle of $\pi$ is a factor of $m$, then there exist
$k_1,k_2,\cdots,k_{d-1}\in [n-1]$ with $d|m$ such that
$k_1,k_2,\cdots,k_{d-1}$ and $n$ form a cycle of $|\pi|$.

If $d=1$, that is $\pi(n)=n^{[j]}$ for some $j$ with $0\leq j\leq
r-1$, then $\pi^m(n)=n^{[jm]}=n$ which implies that $r|jm$. Define
$\pi' \in \mathcal{G}_{r,1,n-1}^m$ by ignoring the last digit of
$\pi$. Then we have
\begin{eqnarray*}
{\rm fix}(\pi)&=&{\rm fix}(\pi')+1,\\
{\rm exc}_A(\pi)&=&{\rm exc}_A(\pi'),\\
\csum(\pi)&=&\csum(\pi')+j.
\end{eqnarray*}

If $d\geq 2$, we can claim that there is $A_{d-1,k}$ cyclic
permutations of length $d$ in $S_d$ with $k$ excedances for $1\leq
k\leq d-1$, where $A_{d-1,k}$ are the Eulerian numbers which are
also the number of permutations on $[d-1]$ with $k-1$ excedances.
This claim can be proved by induction on $d$ and $k$ by noting that
they obey the recurrence relation
\begin{eqnarray}\label{eqn euler}
A_{d,k}=kA_{d-1,k}+(d-k+1)A_{d-1,k-1},
\end{eqnarray}
with the initial values $A_{0,0}=A_{1,1}=1$ and $A_{d,k}=0$ if $d<k$
or $k<1\leq d$. Note that for any cyclic permutation
$C_{d}=(i_1,i_2,\dots,i_d) \in S_{d}$ with $k$ excedances, inserting
$d+1$ into the positions wherever $i_j<i_{j+1}$, there are $k$ ways
to obtain a cyclic permutation $C_{d+1}\in S_{d+1}$ with $k$
excedances, and for any cyclic permutation
$C_{d}=(i_1,i_2,\dots,i_d) \in S_{d}$ with $k-1$ excedances,
inserting $d+1$ into the positions wherever $i_j>i_{j+1}$, there are
$d-k+1$ ways to obtain a cyclic permutation $C_{d+1}\in S_{d+1}$
with $k$ excedances. Conversely, any cyclic permutation $C_{d+1}\in
S_{d+1}$ with $k$ excedances can reduce to the two cases above by
deleting the symbol $d+1$. This analysis makes us get the recurrence
relation (\ref{eqn euler}).

For any cyclic permutation $C$ of length $d$ in $S_d$ with ${\rm
Exc}(C)=\{j\in [d-1]|C(j)>j\}$ such that ${\rm exc}(C)=k$, we can
color the symbols in $C$ with the color set
$\{[0],[1],\cdots,[r-1]\}$ and obtain the colored cyclic permutation
$C'$. Suppose that ${\rm exc}_A(C')=k-i$, we know that ${\rm
Exc}_A(C')\subseteq{\rm Exc}(C)$, which means that ${\rm
exc}(C)-{\rm exc}_A(C)=i$, in other words, there are $i$ number of
symbols in ${\rm Exc}(C)$ with color numbers ranging from $[1]$ to
$[r-1]$, so there are $\binom{k}{i}$ ways to do this.

Let $t=\csum(C')$ and $t_{\ell}$ be the color number of $\ell\in
[d]$, then we have the equation $t=t_1+t_2+\cdots+t_d$ with $0\leq
t_1,t_2,\cdots,t_d\leq r-1$ such that
\begin{itemize}
\item $t_j=0$ for $j\in {\rm Exc}(C)$ and $j$ has a color number
$[0]$, and

\item $1\leq t_j\leq r-1$ for all $j\in {\rm Exc}(C)-{\rm Exc}_A(C')$, so
there are $i$ number of such $j$'s.
\end{itemize}

Therefore there are $U_{d-k,t}^{(i)}$ number of solutions of the
above equation, totally, there are $\binom{k}{i}U_{d-k,t}^{(i)}$
ways to color the symbols in $C$ such that $\csum(C')=t$ and ${\rm
exc}_A(C')=k-i$, where $U_{d-k,t}^{(i)}$ is the coefficient of $x^t$
in $(x+x^2+\cdots+x^{r-1})^i(1+x+\cdots+x^{r-1})^{d-k}$, which can
be expressed as
\begin{eqnarray*}
U_{d-k,t}^{(i)}&=&[x^t](x+x^2+\cdots+x^{r-1})^i(1+x+\cdots+x^{r-1})^{d-k}\\
&=&[x^t]\left(\frac{1-x^r}{1-x}-1\right)^i\left(\frac{1-x^r}{1-x}\right)^{d-k}\\
&=&[x^t]\sum_{j=0}^i(-1)^{i-j}\binom{i}{j}\Big(\frac{1-x^r}{1-x}\Big)^{d+j-k}\\
&=&\sum_{j=0}^i(-1)^{i-j}\binom{i}{j}\sum_{{\ell}\geq
0}(-1)^{{\ell}}\binom{d+j-k}{{\ell}}\binom{d+j+t-k-{\ell}r-1}{t-{\ell}r}.
\end{eqnarray*}
Let $C'=(i_1^{[t_1]},i_2^{[t_2]},\dots,i_d^{[t_d]})$, then
$C'^d=(i_1^{[t]},i_2^{[t]},\dots,i_d^{[t]})$ with
$t=t_1+t_2+\cdots+t_d$, hence
$C'^{m}=(i_1^{[\frac{tm}{d}]},i_2^{[\frac{tm}{d}]},\dots,i_d^{[\frac{tm}{d}]})=1$
implies that $r|\frac{tm}{d}$. For any $\pi \in
\mathcal{G}_{r,1,n}^{m}$ such that the symbol $n$ lies in a cycle
$C'$ of length $d\geq 2$ with $d|m$ (note that there are
$\binom{n-1}{d-1}$ ways to choose the digits of such a cycle),
define $\pi'' \in \mathcal{G}_{r,1,n-d}^{m}$ in the following way:
write $\pi$ in its complete notation, i.e., as a matrix of two rows,
see $(\star)$. The first row of $\pi''$ is $(1,2,\cdots,n-d)$ while
the second row is obtained from the second row of $\pi$ by ignoring
the digits in $C'$ and the other digits are placed with the numbers
$1,2,\dots,n-d$ in an order preserving way with respect to the
second row of $\pi$. The parameters satisfy
\begin{eqnarray*}
{\rm fix}(\pi)&=&{\rm fix}(\pi''),\\
{\rm exc}_A(\pi)&=&{\rm exc}_A(\pi'')+{\rm exc}_A(C'),\\
\csum(\pi)&=&\csum(\pi'')+\csum(C').
\end{eqnarray*}
The above consideration gives the following recurrence
%%\begin{small}
$$\begin{array}{l}
H_{r,1,n}^{(m)}(u,v,w)\\
\qquad=H_{r,1,n-1}^{(m)}(u,v,w)\sum_{\{t|0\leq t<r,
r|tm\}}uw^t+\sum_{d|m,d\geq
2}H_{r,1,n-d}^{(m)}(u,v,w)\binom{n-1}{d-1}A_{m,d}(v,w),
\end{array}$$
%%\end{small}
where
\begin{eqnarray*}
A_{m,d}(v,w)=\sum_{k=1}^{d-1}A_{d-1,k}\sum_{i=0}^k\binom{k}{i}v^{k-i}\sum_{r|\frac{tm}{d}}
U_{d-k,t}^{(i)}w^{t}.
\end{eqnarray*}
Rewriting the recurrence in terms of generating functions, we obtain
that
\begin{eqnarray*}
\lefteqn{\frac{\partial}{\partial x}\mathcal{H}_{r,1}^{(m)}(x;u,v,w)
=\sum_{n \geq1}H_{r,1,n}^{(m)}(u,v,w)\frac{x^{n-1}}{(n-1)!}}\\
&=& \sum_{n \geq
1}\frac{x^{n-1}}{(n-1)!}H_{r,1,n-1}^{(m)}(u,v,w)\sum_{\{t|0\leq t<r,
r|tm\}}uw^t+\\
& & +\sum_{d|m,d\geq
2}A_{m,d}(v,w)\frac{x^{d-1}}{(d-1)!}\sum_{n\geq d}\frac{x^{n-d}}{(n-d)!}H_{r,1,n-d}^{(m)}(u,v,w)\\
&=&\mathcal{H}_{r,1}^{(m)}(x;u,v,w)\Big(\sum_{\{t|0\leq t<r,
r|tm\}}uw^t+\sum_{d|m,d\geq
2}A_{m,d}(v,w)\frac{x^{d-1}}{(d-1)!}\Big).
\end{eqnarray*}

Thus, the generating function $\mathcal{H}_{r,1}^{(m)}(x;u,v,w)$
satisfies
\begin{eqnarray*}
\frac{\frac{\partial}{\partial
x}\mathcal{H}_{r,1}^{(m)}(x;u,v,w)}{\mathcal{H}_{r,1}^{(m)}(x;u,v,w)}
=\sum_{\{t|0\leq t<r, r|tm\}}uw^t+\sum_{d|m,d\geq
2}A_{m,d}(v,w)\frac{x^{d-1}}{(d-1)!}.
\end{eqnarray*}
Integrating with respect to $x$ on both sides of the above
differential equation, using the fact that
$\mathcal{H}_{r,1}^{(m)}(0;u,v,w)=1$, we obtain the explicit
expression for $\mathcal{H}_{r,1}^{(m)}(x;u,v,w)$ given in Theorem
\ref{maintheo}, and hence we complete the proof.\qed\vskip0.2cm

Specially, if $m=p$ is a prime, then we have
\begin{corollary}\label{coro3.2}
Let $r\geq 1$ and $p$ be a prime. The generating function
$\mathcal{H}_{r,1}^{(p)}(x;u,v,w)$ is given by
\begin{eqnarray*}
\exp\left\{ux\lambda_{r,p}(w)
+\frac{x^p}{p!}\sum_{k=1}^{p-1}A_{p-1,k}\sum_{i=0}^k\binom{k}{i}v^{k-i}\sum_{j\geq
0}U_{p-k,jr}^{(i)}w^{jr}\right\},
\end{eqnarray*}
where $A_{p-1,k}$ is the Eulerian number, $U_{p-k,jr}^{(i)}$ is the
coefficient of $x^{jr}$ in
$$(x+x^2+\cdots+x^{r-1})^i(1+x+\cdots+x^{r-1})^{p-k},$$
$\lambda_{r,p}(w)=\sum_{i=0}^{p-1}{w^{\frac{ir}{p}}}$ for $p|r$, and
$\lambda_{r,p}(w)=1$ for $p\not| r$.
\end{corollary}

For the sake of comparison, the cases $p=2$ and $p=3$ in Corollary
\ref{coro3.2} generate the explicit formulas for
$\mathcal{H}_{r,1}^{(2)}(x;u,v,w)$ and
$\mathcal{H}_{r,1}^{(3)}(x;u,v,w)$, that is
\begin{eqnarray*}
\mathcal{H}_{r,1}^{(2)}(x;u,v,w)&=& \exp({ux\lambda_{r,2}(w)
+\frac{x^2}{2}(v+(r-1)w^{r})}),\\
\mathcal{H}_{r,1}^{(3)}(x;u,v,w)&=& \exp({ux\lambda_{r,3}(w)
+\frac{x^3}{6}B_{3,3}(v,w)}),
\end{eqnarray*}
where
$B_{3,3}(v,w)=v^2+v(1+3(r-1)w^r)+(r^2-1)w^{r}+(r-1)(r-2)w^{2r}$.
%%and $A_{4,2}(v,w)=v+(r-1)w^r$,
%%$A_{4,4}(v,w)=v^3+v^2(4+7(r-1)w^r)+v(1+2(r-1)(r+2)w^r+5(r-1)(r-2)w^{2r})+(r^3-1)w^{r}+2(r^2-1)(2r-3)w^{2r}+2\binom{r}{3}w^{3r}$
%%for $2\not|r$

Now let us compute the exponential generating function
$\mathcal{H}_{r,s}^{(m)}(x;u,v,w)$ for the sequence
$\{H_{r,s,n}^{(m)}(u,v,w)\}_{n\geq0}$. For any $\sigma\in
\mathcal{G}_{r,s,n}^{m}$, we have $\csum(\sigma)\equiv 0  \; ({\rm
mod} \; s)$, so we should collect all the terms in which the
exponent of $w$ in $\mathcal{H}_{r,1}^{(m)}(u,v,w)$ is a
multiplication of $s$. This observation can make us get the
following
\begin{theorem}\label{theo3.4}
Let $r,m,s\geq 1$, define
${\mathcal{H}_{r,1}^{(m)}(x;u,v,yw)}=\sum_{n\geq
0}G_{m,r,n}(x;u,v,w)y^n.$ Then
\begin{eqnarray*}
\mathcal{H}_{r,s}^{(m)}(x;u,v,w)&=&\sum_{k\geq
0}G_{m,r,sk}(x;u,v,w).
\end{eqnarray*}
\end{theorem}

Now let us focus on the case $m=2$. Recall that
\begin{eqnarray*}
{\mathcal{H}_{r,1}^{(2)}(x;u,v,w)}=\left\{
\begin{array}{ll}
e^{ux+\frac{1}{2}x^2(v+(r-1)w^{r})},& {\rm if}\ r\ {\rm odd},\\
e^{ux(1+w^{\frac{r}{2}})+\frac{1}{2}x^2(v+(r-1)w^{r})},& {\rm if}\
r\ {\rm even\ }.
\end{array}\right.
\end{eqnarray*}
Then by Theorem \ref{theo3.4}, we can compute the explicit formula
for $\mathcal{H}_{r,s}^{(2)}(x;u,v,w)$. Since $s\mid r$, we have two
cases either $r$ odd or $r$ even.
\begin{itemize}
\item If $r$ is an odd number, then it is clear that the exponent of $y$ in each term
of the expansions of $\mathcal{H}_{r,1}^{(2)}(x;u,v,yw)$ is always a
multiplication of $s$. Hence,
$$\mathcal{H}_{r,s}^{(2)}(x;u,v,w)=\mathcal{H}_{r,1}^{(2)}(x;u,v,w).$$

\item Similarly, if $r$ is an even number and $s|\frac{r}{2}$, we have that
$$\mathcal{H}_{r,s}^{(2)}(x;u,v,w)=\mathcal{H}_{r,1}^{(2)}(x;u,v,w).$$

\item Let $r$ be any even number such that $s\nmid\frac{r}{2}$. Since
$e^{ux(1+(yw)^{\frac{r}{2}})}=e^{ux}\sum_{k\geq0}\frac{(ux(yw)^{\frac{r}{2}})^k}{k!}$
and
$e^{\frac{1}{2}x^2(v+(r-1)(yw)^{r})}=e^{\frac{1}{2}x^2v}\sum_{k\geq0}\frac{((r-1)x^2(yw)^r)^k}{2^kk!}$,
then by collecting the coefficients of $y$ in
$\mathcal{H}_{r,1}^{(2)}(x;u,v,w)$ such that the exponent $y$ is a
multiplication of $s$, we get that
$$e^{\frac{1}{2}x^2(v+(r-1)(yw)^{r})}\sum_{k\geq
0}\frac{(ux)^{2k}(yw)^{kr}}{(2k)!}=e^{ux+\frac{1}{2}x^2(v+(r-1)(yw)^{r})}\frac{e^{uxw^{\frac{r}{2}}}+e^{-uxw^{\frac{r}{2}}}}{2}.$$
\end{itemize}
Therefore, the above cases gives the following result.
\begin{proposition}\label{pros1}
We have
\begin{eqnarray*}
{\mathcal{H}_{r,s}^{(2)}(x;u,v,w)}&=&\left\{
\begin{array}{ll}
e^{ux+\frac{1}{2}x^2(v+(r-1)w^{r})},& {\rm if}\ r\ {\rm odd},\\
e^{ux(1+w^{\frac{r}{2}})+\frac{1}{2}x^2(v+(r-1)w^{r})},& {\rm if}\
r\ {\rm even\ and\ }
s\nmid\frac{r}{2},\\
e^{ux+\frac{1}{2}x^2(v+(r-1)w^{r})}\frac{e^{uxw^{\frac{r}{2}}}+e^{-uxw^{\frac{r}{2}}}}{2},&
{\rm if}\ r\ {\rm even\ and\ } s\nmid \frac{r}{2}.
\end{array}\right.
\end{eqnarray*}
\end{proposition}

Note that $\mathcal{H}_{r,s}^{(2)}(x;u,v,w)$ is the generating
function for the number of involutions in
$\mathcal{G}_{r,s,n}^{(2)}$. By expanding the generating functions,
Bagno, Garber and Mansour \cite{BGM} obtained the explicit formulas
for the number of involutions in $\mathcal{G}_{r,s,n}^{(2)}$. But
the expression in Proposition 5.7 \cite{BGM} should be corrected by
the third case of $\mathcal{H}_{r,s}^{(2)}(x;u,v,w)$ and hence
Corollary 5.8, 5.9 and 5.10 therein should be the following three
corollaries, respectively.
\begin{corollary}
The polynomial $H_{r,s,n}^{(2)}(u,v,w)$ is given by
\begin{eqnarray*}
\sum_{k_1+2k_2+2k_3=n}\frac{n!}{k_1!(2k_2)!k_3!}\cdot\frac{u^{k_1+2k_2}w^{rk_2}(v+(r-1)w^r)^{k_3}}{2^{k_3}}.
\end{eqnarray*}
\end{corollary}
\begin{corollary} Let $r\geq1$.  The number of colored involutions in
$\mathcal{G}_{r,s,n}^{(2)}$ ($r$ is even, $s \nmid \frac{r}{2}$)
with exactly $k$ absolute fixed points and ${\rm exc}_A (\pi)=\ell$
is given by
\begin{eqnarray*}
\sum_{k+2k_3=n,k_1+2k_2=k}\binom{k_3}{\ell}\cdot\frac{n!}{k_1!(2k_2)!k_3!}\cdot
\frac{(r-1)^{k_3-\ell}}{2^{k_3}}.
\end{eqnarray*}
\end{corollary}
\begin{corollary}
The number of involutions  $\pi \in \mathcal{G}_{r,s,n}^{(2)}$ ($r$
is even, $s \nmid \frac{r}{2}$) with ${\rm exc}^{{\rm Clr}}(\pi)=k$
is given by
\begin{eqnarray*}
\sum_{k_1+2k_2+2k_3=n,\
r(k_2+k_3)=k}\frac{n!}{k_1!(2k_2)!k_3!}\cdot\left(\frac{r}{2}\right)^{\frac{k}{r}}.
\end{eqnarray*}
\end{corollary}

{\bf Acknowledgment} The authors would like to thank Eli Bagno and
David Garber for reading previous version of the present paper and
for a number of helpful discussions.
%==============================================================================================================


\begin{thebibliography}{99}

\bibitem{BG}
E. Bagno and D. Garber, {\it On the excedance number of colored
permutation groups}, Semi. Loth. Comb. {\bf 53} (2006), Art. B53f,
17 pp. (Electronic, available at http://igd.univ-lyon1.fr/\~{}slc).

\bibitem{BGM}
E. Bagno, D. Garber and T. Mansour, {\it Excedance number for
involutions in complex reflection groups}, Semi. Loth. Comb. {\bf
56} (2007), Art. B56d, 11 pp. (Electronic, available at
http://igd.univ-lyon1.fr/\~{}slc).


\bibitem{KZ}
G.\ Ksavrelof and J.\ Zeng, {\it Two involutions for signed
excedance numbers}, Semi. Loth. Comb. {\bf 49} (2003), Art. B49e, 8
pp. (Electronic, available at http://igd.univ-lyon1.fr/\~{}slc).

%\bibitem{ST} G, C. Shephard and J. A. Todd, \emph{Finite unitary reflection groups},
%Canad. J. Math., \textbf{6} (1954), 274-304.

\bibitem{SS}
E. Steingr{\'{\i}}msson, {\it Permutation statistics on indexed
permutations}, Europ. J. Comb. {\bf 15:2} (1994) 187--205.
\end{thebibliography}
\end{document}